\documentclass[12pt]{article}
\usepackage{amsfonts}
\usepackage{amsmath}
\usepackage{amssymb}
\usepackage{oldlfont}
\newtheorem{Theorem}{Theorem}[section]
\newtheorem{Lemma}[Theorem]{Lemma}
\newtheorem{Proposition}[Theorem]{Proposition}

\newtheorem{Corollary}[Theorem]{Corollary}

\begin{document}
\title{\bf{Chaining Techniques and their Application to Stochastic Flows}}
\author{Michael Scheutzow}


\maketitle
\begin{abstract}
We review several competing {\em chaining} methods to estimate the supremum, the diameter of the range or the 
modulus of continuity of a stochastic process in terms of tail bounds of their two-dimensional 
distributions. 
Then we show how they can be applied to obtain upper bounds for the growth of 
bounded sets under the action of a stochastic flow. 
\end{abstract}

\section{Introduction}
Upper and lower bounds for the (linear) growth rates of the 
diameter of the image of a bounded set in ${\mathbf R}^d$ 
under the action of a stochastic flow under various conditions have 
been shown in \cite{cs,css1,css2,ls01,ls03,ss}. 
In this survey, we will only discuss upper bounds. A well-established class of methods to obtain probability bounds 
for the supremum of a process are {\em chaining} techniques. Typically they transform bounds for the one- and two-dimensional 
distributions of the process into upper bounds of the supremum (for a real-valued process) 
 or the diameter of the range of the process  (for a process taking values in a metric space). In the next section, we will 
present some of these techniques, the best-known being Kolmogorov's continuity theorem, which not only states the 
existence of a continuous modification, but also provides explicit probabilistic upper bounds for the modulus of 
continuity and the diameter of the range of the process. We will also state a result which we call {\em basic chaining}. 
Further we will briefly review some of the results from Ledoux and Talagrand \cite{lt} and a rather general version 
of the {\em GRR-Lemma} named after Garsia, Rodemich and Rumsey \cite{grr}. Except for the result of Ledoux and Talagrand, 
we will provide proofs for the chaining lemmas in the appendix (in order to keep the article reasonably self-contained but 
also because we chose to formulate the chaining results slightly differently compared to the literature). We wish to 
point out however, that nothing in that section is essentially new and that it is not meant to be a complete survey 
about chaining. The reader who is interested in learning more about chaining should consult the literature,  for 
example the monograph by Talagrand \cite{T05}.  

In order to obtain good upper bounds on the diameter of the image of a bounded set ${\mathcal X}$ under a 
stochastic flow $\phi$ which 
is generated (say) by a stochastic differential equation on ${\mathbf R}^d$ with coefficients which are bounded and 
Lipschitz continuous, one can try to apply the chaining techniques directly to the process $\phi_{0,T}(x)$, 
$x \in {\mathcal X}$. 
This is what we did in \cite{css2} using basic chaining. It worked, but it was a nightmare (for the reader, the referee and 
us). The reason was, that the two-point motion of such a flow behaves quite differently depending on 
whether the two points are very close (then the Lipschitz constants determine the dynamics) or not (then the bounds on 
the coefficients do). This requires a rather sophisticated choice of the parameters or functions in the chaining lemmas. 
The papers \cite{ls01,ls03} provided somewhat simpler proofs using the chaining methods of Ledoux-Talagrand and the GRR-Lemma 
respectively. The approach presented here is (in our opinion) much simpler and transparent than the previous ones. The reason 
is, that we strictly separate the local and the global behaviour in the following sense: for a given (large) time $T$ and 
a positive number $\gamma$, we cover the set ${\mathcal X}$ with balls (or cubes) of radius $\exp\{-\gamma T\}$. For each  
center of such a ball, we estimate the probability, that it leaves a ball with radius $\kappa T$ around zero up to time $T$
using large deviations estimates (Proposition \ref{1pt}). This probability bound depends only on the bounds on the 
coefficients and not on the Lipschitz constants. In addition, we provide an upper bound for the probablity that a 
particular one  of the small balls achieves a diameter of 1 (or some other fixed positive number) up to time $T$ 
(Theorem \ref{SB}).  
This bound only involves the Lipschitz constants and not the bounds on the coefficients. To obtain such a bound, we use 
chaining. We will allow ourselves the luxury of five proofs of this result using each of the chaining methods -- 
with even two proofs using the GRR-Lemma. Since we are only interested in the behaviour of the image of a very small 
ball up to the time its radius becomes 1, things become much easier compared to the approach in 
\cite{css2} mentioned above. In fact, we can use a polynomial function $\Psi$ when applying 
the GRR-Lemma or the LT-Lemma and this is why Kolmogorov's Theorem, which also uses polynomial moment bounds, 
turns out to be just as efficient as the other (more sophisticated) methods.  The proof of Theorem \ref{diam}, which 
provides an explicit upper bound for the linear growth rate, now becomes almost straightforward: the probability that the 
diameter of the image of ${\mathcal X}$ under the flow up to time $T$ exceeds $\kappa T$ is bounded from above by the number of 
small balls multiplied by the (maximal) probability that a center reaches a modulus of $\kappa T-1$ or the diameter 
of a small ball exceeds 1. This bound -- which is still a function of the parameter $\gamma$ -- turns out to be exponentially 
small in $T$ provided $\kappa$ is large enough and $\gamma$ is chosen appropriately. An application of the first 
Borel-Cantelli Lemma then completes the proof.

We talk about stochastic flows above, but the results are true under less restrictive conditions. For the upper bound 
of the growth of a small ball (Theorem \ref{SB}), it suffices that
the underlying motion $\phi_t(x)$ is jointly continuous and that (roughly speaking) 
the distance of two trajectories does not grow 
faster than a geometric Brownian motion (this is hypothesis (H) in Section \ref{work}). 
In the special case of a (spatially) differentiable and translation invariant Brownian flow, 
Theorem  \ref{SB} can be improved slightly. This is shown in Theorem \ref{diff}. Its proof is completely different 
from that of Theorem  \ref{SB}: it does not use any chaining whatsoever.

\section{The Competitors}
In the following, we will always assume that $(\hat E,\hat \rho)$ is a complete, separable, metric space. 
Further, $\Psi: [0,\infty) \to [0,\infty)$ will always be a strictly increasing function which satisfies 
$\Psi(0)=0$. If -- in addition -- $\Psi$ is convex, then it is called a {\em Young function}. 
For a Young function $\Psi$, one defines the corresponding {\em Orlicz norm} of a
real-valued random variable $Z$ by 
$$
\|Z\|_{\Psi}:=\inf\{c>0: {\mathbf E} \Psi(|Z|/c) \le 1\}.
$$
We will also need a totally bounded metric space $(\Theta ,d)$ with diameter $D>0$. The minimal number of closed balls of 
radius $\varepsilon$ needed to cover $\Theta$ will be denoted by $N(\Theta,d;\varepsilon)$ and will be called 
{\em covering numbers}.
A finite subset $\Theta_0$ of $\Theta$ is called an {\em $\varepsilon$-net}, if $d(x,\Theta_0) \le \varepsilon$ for each 
$x \in \Theta$ (we use $x$ rather than $t$, because in our application $\Theta$ will be a subset of the space ${\mathbf R}^d$). 
We will abbreviate
$$
J:=\int_0^D \Psi^{-1}(N(\Theta,d;\varepsilon)) \,{\mathrm d}\varepsilon.
$$
Further, let $Z_x$, $x \in \Theta$ be an $\hat E$-valued process on some probability space $(\Omega,{\cal {F}},{\mathbf P})$.
We will denote the Euclidean norm, the $l_1$-norm and the maximum norm on ${\mathbf R}^d$ by $|.|$, $|.|_1$ and $|.|_{\infty}$ 
respectively. Whenever a constant is denoted by $c$ with some index, then its value can change from line to line. 
We start with the well-known continuity theorem of Kolmogorov.

\begin{Lemma} {\bf (Kolmogorov)} \label{Kolmo} 
Let $\Theta=[0,1]^d$ and assume that there exist $a,b,c >0$ such that for all $x,y \in [0,1]^d$, we have
$$
{\mathbf E} \left((\hat \rho(Z_x,Z_y))^a\right) \le c |x-y|_1^{d+b}.
$$
Then $Z$ has a continuous modification (which we denote by the same symbol). For each $\kappa \in (0,b/a)$, there exists a 
random variable $S$ such that ${\mathbf E}(S^a) \le \frac{cd 2^{a\kappa-b}}{1- 2^{a\kappa-b}}$ and
$$
\sup \left\{ \hat \rho(Z_x(\omega),Z_y(\omega)): x,y \in [0,1]^d,|x-y|_{\infty} \le r \right\} \le 
\frac{2d}{1-2^{-\kappa}} S(\omega) r^{\kappa}
$$
for each $r \in [0,1]$.
In particular, for all $u > 0$, we have
\begin{equation}\label{estikolmo}
{\mathbf P}\left\{ \sup_{x,y \in [0,1]^d} \hat \rho(Z_x,Z_y) \ge u \right\} \le 
\left( \frac{2d}{1-2^{-\kappa}}\right)^a \frac{cd 2^{a\kappa-b}}{1- 2^{a\kappa-b}}u^{-a}.
\end{equation}
\end{Lemma}

\begin{Lemma} {\bf (Basic Chaining)} \label{Chaining}
Let $Z$ have continuous paths. 
Further, let $\delta_j$, $j=0,1,...$ be a sequence of positive real numbers such that $\sum_{j=0}^{\infty} \delta_j <\infty$ and let 
$\Theta_j$ be a $\delta_j$-net in $(\Theta,d)$, $j=0,1,2,...$ such that $\Theta_0=\{x_0\}$ is a singleton.

Then for any $u>0$ and any sequence of positive $\varepsilon_j$ with $\sum_{j=0}^{\infty} \varepsilon_j\le 1$,
$$
{\mathbf P} \left\{ \sup_{x,y \in \Theta} \hat \rho (Z_x,Z_y) \ge u \right\} 
\le \sum_{j=0}^{\infty} |\Theta_j| \sup_{d(x, y)   \le \delta_j}  {\mathbf P}\left\{\hat \rho (Z_x,Z_y) \ge \varepsilon_j u/2  \right\}.
$$
\end{Lemma}

The following lemma combines Theorems 11.1., 11.2., 11.6., and (11.3) in \cite{lt} (observe the obvious typo in
(11.3) of \cite{lt}: $\psi^{-1}$ should be replaced by $\psi$).
\begin{Lemma} {\bf (LT-Chaining)} \label{LT}
Let $\Psi$ be a Young function such that $J < \infty$. Assume that there exists a constant $c>0$ such that for 
all $x,y \in \Theta$
$$
\|\hat \rho(Z_x,Z_y )\|_{\Psi} \le c d(x,y).
$$
Then $Z$ has a continuous modification (which we denote by the same symbol). Further, for each set $A \in {\cal {F}}$, we have
$$
\int_A \sup_{x,y \in \Theta} \hat \rho(Z_x,Z_y)  \,{\mathrm d} {\mathbf P} \le 8 {\mathbf P} (A) c \int_0^D  
\Psi^{-1} \left( \frac{N(\Theta,d;\varepsilon)}{{\mathbf P} (A)}
\right)  {\mathrm d} \varepsilon.
$$ 
If, in addition, there exists $c_{\Psi} \ge 0$ which satisfies
$\Psi^{-1}(\alpha \beta) \le c_{\Psi} \Psi^{-1}(\alpha) \Psi^{-1}(\beta)$ for all $\alpha,\beta \ge 1$, then for all $u > 0$, we have
$$
{\mathbf P}\left\{ \sup_{x,y \in \Theta} \hat \rho(Z_x,Z_y) \ge u \right\} \le \left( \Psi\left( \frac{u}{8c c_{\Psi} J} \right)\right)^{-1}. 
$$
\end{Lemma}

The following version of the GRR-Lemma seems to be new. It is a joint upgrade (up to  constants) 
of  \cite{dz}, Theorem B.1.1 and \cite{ai}, Theorem 1. Even though the version in  \cite{dz} meets  our 
demands, we present a more general version below and prove it in the appendix.

\begin{Lemma} {\bf (GRR)} \label{GRR}
Let $(\Theta,d)$ be an arbitrary metric space (not necessarily totally bounded), $m$ a measure on the Borel sets of $\Theta$ 
which is finite on bounded subsets and
let $p:[0,\infty) \to [0,\infty)$ be continuous and strictly increasing and $p(0)=0$. If
$f:\Theta \to \hat E$ is continuous such that
$$
V:=\int_{\Theta}  \int_{\Theta} \Psi \left( \frac{\hat \rho(f(x),f(y))}{p(d(x,y))} \right)  \,{\mathrm d} m(x)   \,{\mathrm d}m(y) < \infty,
$$
then we have
\begin{itemize}
\item[\textrm{(i)}] $\hat \rho(f(x),f(y)) \le 8 \max_{z \in \{x,y\}} 
\int_0^{4d(x,y)} \Psi^{-1} \left( \frac{4V}{m(K_{s/2}(z))^2} \right)  {\mathrm d} p(s)$,
\item[\textrm{(ii)}] $\hat \rho(f(x),f(y)) \le 8 N \max_{z \in \{x,y\}} 
\int_0^{4d(x,y)} \Psi^{-1} \left( \frac{4}{m(K_{s/2}(z))^2} \right)  {\mathrm d}p(s)$,
\end{itemize}
where 
$$
N:=\inf\{\kappa>0: \int \int \Psi \left( \frac{1}{\kappa} \frac{\hat \rho(f(x),f(y))}{p(d(x,y))}\right) 
{\mathrm d} m(x)   \,{\mathrm d}m(y) \le 1\}
$$
and $K_s(z)$ denotes the closed ball with center $z$ and radius $s$.
In the definition of $V$ and $N$, $0/0$ is interpreted as zero, while in the conclusions $V/0$ and $N \times \infty$ are 
interpreted as $\infty$ even if $V=0$ or $N=0$.
\end{Lemma}

\noindent{\bf Remark} \hspace{.3cm} When applying one of the chaining methods above, one is forced to 
choose the function $\Psi$ 
(for LT-chaining and GRR) or other parameters (in basic chaining and Kolmogorov's Theorem). One might suspect that 
it is wise to choose $\Psi$ in such a way, that it increases as quickly as possible subject to the constraint that $J < \infty$ 
(in LT-chaining) because this will guarantee sharper tail estimates for the suprema in question. It may therefore come as a
surprise that we will be able to obtain optimal estimates by choosing polynomial functions $\Psi$ and that Kolmogorov's 
Theorem, which only allows for polynomial functions, will be just as good as the much more sophisticated LT-chaining (for example).
The reason for this is, that we will use chaining only to estimate the probability that the diameter of the image of a small ball 
under a flow (for example) exceeds a fixed value (for example 1) up to a given time $T$ and we do not care how large the diameter 
is if it exceeds this value.\\

\noindent{\bf Remarks about the chaining literature.} The GRR-Lemma was first published in \cite{grr} in the special 
case $\Theta=[0,1]$. A version where $\Theta$ is an open bounded set in ${\mathbf R}^d$ can be found in \cite{dz}, 
Appendix B (with $m=$ Lebesgue measure). Walsh (\cite{walsh}, Theorem 1.1) requires $\Theta=[0,1]^d$, 
$m=$ Lebesgue measure, $\Psi$ convex and $f$ real-valued but does not assume that $f$ be continuous. 
The GRR Lemma in \cite{ai} is similar to ours but they assume that $p=$ identity. 
Dalang et.~al.~\cite{dkn} prove a version which is also similar to ours. They assume that the function $\Psi$ is 
convex (which we don't) and in turn obtain a smaller multiplicative constant. Like Walsh \cite{walsh}, they do not 
need to assume that the function $f$ is continuous.

Lemma \ref{Chaining} appeared (in a slightly different form) in \cite{css2}, but even at that time it was adequate to 
call it {\em essentially well-known}. Indeed, the idea of choosing a sequence of finite $\delta$-nets with 
$\delta \to 0$ is at the heart of the chaining method (see, e.g.~\cite{p}).

One can find more general {\em anisotropic} versions of Kolmogorov's continuity Theorem \ref{Kolmo} in which 
the right hand side $c|x-y|_1^{d+b}$ is replaced by $c\sum_{i=1}^d|x_i-y_i|^{\alpha_i}$ where $\sum_{i=1}^d 
\alpha_1^{-1} <1$, see, e.g.~\cite{k} or  \cite{dkn}. We point out that Kolmogorov's Theorem can be regarded 
as a corollary (possibly up to multiplicative constants) of both LT-Chaining (Lemma \ref{LT}) and certain variants 
of the GRR Lemma, see \cite{lt} and \cite{walsh} respectively.

\section{Chaining at Work}\label{work}

Let $(t,x) \mapsto \phi_t(x)$ be a continuous random field, $(t,x) \in [0,\infty) \times {\mathbf R}^d$ 
taking values in a separable complete metric space $(E,\rho)$. We will always assume 
that $\phi$ satisfies the following condition:\\

\noindent {\bf (H)}: There exist $\Lambda \ge 0$, $\sigma >0$ and $\bar c>0$ such that for each $x,y \in {\mathbf R}^d$,
$T>0$, and $q \ge 1$, we have
$$
\left( {\mathbf E} \sup_{0 \le t \le T} (\rho(\phi_t(x),\phi_t(y)))^q \right)^{1/q} \le \bar c \,|x-y| \, \exp\{(\Lambda + \frac{1}{2}q \sigma^2)T\}.
$$

A sufficient condition for (H) to hold (with $\bar c = 2$) is the following condition (H').\\

\noindent {\bf (H')}: There exist $\Lambda \ge 0$, $\sigma >0$ such that for each $x,y \in {\mathbf R}^d$, 
there exists a standard Brownian motion $W$, such that
\begin{equation}\label{H}
\rho(\phi_t(x),\phi_t(y)) \le |x-y|\, \exp \{ \Lambda t + \sigma W_t^* \},
\end{equation}
where $W_t^* :=\sup_{0 \le s \le t} W(s)$.\\

We will verify in Lemma \ref{sde} that (H') and hence (H) is satisfied
for the solution flow of a stochastic differential equation on ${\mathbf R}^d$ with global Lipschitz 
coefficients.

If there exists some $\nu>0$ such that (\ref{H}) holds only for $t \le \inf\{s \ge 0: \rho(\phi_s(x),\phi_s(y)) \ge \nu\}$, then 
(H') holds provided that $\rho$ is replaced by the metric $\bar \rho(x_1,x_2):=\rho(x_1,x_2)\wedge \nu$. Choosing $\nu$ small 
allows in some cases to use smaller values of $\Lambda$ and/or $\sigma$ and thus to improve the asymptotic bounds in the 
following theorem.

In fact the application of Lemma \ref{Kolmo} or \ref{LT} below shows that the existence of a continuous modification of $\phi$ 
w.r.t. $x$ follows from (H).

In the following Theorem, we will provide an upper bound for the probability that the image of 
a ball which is exponentially small in $T$, attains diameter 1 (say) up to time $T$. 

\begin{Theorem}\label{SB} Assume $\mathrm {(H)}$ and let $\gamma > 0$. Define
\begin{eqnarray*}
I(\gamma) :=   \left\{ 
\begin{array}{ll}
\frac{(\gamma-\Lambda)^2}{2\sigma ^2} &{\rm if } \, \gamma \ge \Lambda + \sigma^2 d\\
d (\gamma-\Lambda - \frac{1}{2}\sigma^2 d)  &{\rm if } \,\Lambda + \frac{1}{2}\sigma^2 d \le \gamma \le \Lambda + \sigma^2 d\\
0 &{\rm if } \,\gamma \le \Lambda + \frac{1}{2}\sigma^2 d.
\end{array}
 \right.
\end{eqnarray*}
Then, for each $u>0$, we have
$$
\limsup_{T \to \infty} \frac{1}{T} \sup_{{\mathcal X}_T} \log {\mathbf P}\{\sup_{x,y \in {\mathcal X}_T} \sup_{0 \le t \le T} \rho(\phi_t(x),\phi_t(y)) \ge u\} \le -I(\gamma),
$$
where $\sup_{{\mathcal X}_T}$ means that we take the supremum over all cubes ${\mathcal X}_T$ in ${\mathbf R}^d$ 
with side length $\exp\{-\gamma T\}$.
\end{Theorem}

We will first provide five different proofs of Theorem \ref{SB} by using Lemmas \ref{Kolmo}, 
\ref{Chaining}, \ref{LT}, and \ref{GRR} respectively. We will always use the space 
$\hat E=C([0,T],E)$ equipped with the sup-norm $\hat \rho$, where $(E,\rho)$ is a complete separable metric space as above.\\

\noindent{\bf Proof of Theorem \ref{SB} using Lemma \ref{Kolmo}.} Let $T>0$. 
Without loss of generality, we assume that ${\mathcal X}:={\mathcal X}_T=[0,{\mathrm e}^{-\gamma T}]^d$. 
Define $Z_x(t):=\phi_t({\mathrm e}^{-\gamma T} x)$, $x \in {\mathbf R}^d$. 
For $q \ge 1$, (H) implies
$$
\left({\mathbf E} \sup_{0 \le t \le T} \rho(Z_x(t),Z_y(t))^q\right)^{1/q} \le \bar c {\mathrm e}^{-\gamma T} |x-y|  
{\mathrm e}^{(\Lambda+\frac{1}{2}q\sigma^2)T},
$$ 
i.e.~ the assumptions of Lemma \ref{Kolmo} are satisfied with $a=q$, $c=\bar c^q \exp\{(\Lambda-\gamma + \frac{1}{2} q \sigma^2)qT\}$ 
and $b=q-d$ for any $q>d$. Therefore we get for $\kappa \in (0,b/a)$:
\begin{align*}
&{\mathbf P}\{\sup_{x,y \in {\mathcal X}} \sup_{0 \le t \le T} \rho(\phi_t(x),\phi_t(y))\ge u\} \\ 
&\hspace{3cm} \le \left( \frac{2d}{1-2^{-\kappa}}\right)^q 
\frac{\bar c^q d 2^{a\kappa - b}}{1- 2^{a\kappa - b}} \exp\{(\Lambda-\gamma+\frac{1}{2} q \sigma^2)qT\} u^{-q}.
\end{align*}
Taking logs, dividing by $T$, letting $T \to \infty$ and optimizing over $q>d$ yields Theorem \ref{SB}. \hfill $\Box$\\

\noindent{\bf Proof of Theorem \ref{SB} using Lemma \ref{Chaining}.}
Let $\gamma>\Lambda$ and take a cube $\Theta={\mathcal X}_T$ of side length $\exp\{-\gamma T\}$. 
Then we apply the Chaining Lemma \ref{Chaining} to $\Theta$ with $\delta_j=\exp\{-\gamma T\} \sqrt{d} 2^{-j-1}$ 
and $\varepsilon_j=C/(j+1)^2$, $j=0,1,...$, where the constant $C$ is chosen such that the $\varepsilon_j$ sum up to 1. 
Then there exist subsets $\Theta_j$ of $B$ with cardinality $|\Theta_j|=2^{jd}$ such that
the assumptions of Lemma \ref{Chaining} are satisfied. In particular, $x_0$ is the center of the cube $\Theta$. 
For $q>d$, we get

\begin{eqnarray}
&&{\mathbf P}\left\{ \sup_{x,y \in \Theta} \sup_{t \in [0,T]} \rho ( \phi_t(x),\phi_t(y) ) \ge u \right\}\nonumber\\
&\le&\sum_{j=0}^{\infty} 2^{dj} \sup_{|x-y| \le \delta_j} {\mathbf P} \left\{ \sup_{t \in [0,T]} \rho( \phi_t(x),\phi_t(y)) 
\ge \varepsilon_j u/2 \right\}\nonumber\\ 
&\le&\sum_{j=0}^{\infty} 2^{dj} (\varepsilon_j u/2)^{-q} 
\sup_{|x-y| \le \delta_j} {\mathbf E}\left(\sup_{0 \le t \le T} \rho(\phi_t(x),\phi_t(y)) \right)^q\nonumber\\
&\le&{\mathrm e}^{((\Lambda-\gamma)q+\frac{1}{2}\sigma^2 q^2)T} \bar c^q d^{q/2} u^{-q}  
\sum_{j=0}^{\infty} 2^{(d-q)j}\varepsilon_j^{-q}.\label{estimate}
\end{eqnarray}
The sum converges since $q>d$ and the $\varepsilon_j$ decay polynomially. Taking logs in (\ref{estimate}), 
dividing by $T$, letting $T \to \infty$ and optimizing over $q>d$ yields Theorem \ref{SB}. \hfill $\Box$\\

\noindent{\bf Proof of Theorem \ref{SB} using Lemma \ref{LT}.}
Fix $T>0$ and $q>d$. We apply Lemma \ref{LT} with $\Psi(x)=x^q$ (then $c_{\Psi}=1$). 
Inequality (H) shows that the assumptions are 
satisfied with $c=\bar c \exp\{(\Lambda + \frac{1}{2} q\sigma^2)T\}$. Further, we have
$$
J:=\int_0^{\sqrt{d} {\mathrm e}^{-\gamma T}} N([0,{\mathrm e}^{-\gamma T}]^d,|.|;\varepsilon)^{1/q} \,{\mathrm d}  
\varepsilon \le c_{d,q} {\mathrm e}^{-\gamma T}.
$$
Therefore, we obtain
\begin{eqnarray*}
{\mathbf P}\{\sup_{x,y \in {\mathcal X}} \sup_{0 \le t \le T} \rho(\phi_t(x),\phi_t(y))\ge u\} \le (8 J c)^q u^{-q} \hspace{1cm}\\
\hspace{1cm}\le \bar c^q \tilde c_{d,q} \exp\{(\Lambda q -\gamma q +\frac{1}{2} q^2 \sigma^2)T\}  u^{-q}
\end{eqnarray*}
Taking logarithms, dividing by $T$, letting $T \to \infty$ and optimizing over $q>d$ yields the claim in Theorem \ref{SB}.
\hfill $\Box$\\

\noindent{\bf Proof of Theorem \ref{SB} using Lemma \ref{GRR}.}
Let $p(s):=s^{(2d+\varepsilon)/q}$, where $\varepsilon \in (0,1)$ and $q > d +\varepsilon$. Define
$$
V:=\int_{{\mathcal X}_T} \int_{{\mathcal X}_T} \sup_{0 \le t \le T} \frac{\rho(\phi_t(x),\phi_t(y))^q}{p(|y-x|)^q} 
 \,{\mathrm d}  x  \,{\mathrm d}  y.
$$
Let $m$ be Lebesgue measure restricted to ${\mathcal X}_T$ and $\Psi(x)=x^q$. By (H),
\begin{eqnarray*}
{\mathbf E} V &\le& \bar c^q {\mathrm e}^{(\Lambda +\frac{1}{2}\sigma^2 q) q T} 
\int_{{\mathcal X}_T} \int_{{\mathcal X}_T} |y-x|^{q-2d-\varepsilon}  \,{\mathrm d}  x  \,{\mathrm d}  y\\
 &\le& \bar c^q   {\mathrm e}^{(\Lambda +\frac{1}{2}\sigma^2 q) q T} 
{\mathrm e}^{-\gamma d T} \int_{\{|y| \le \sqrt{d}{\mathrm e}^{-\gamma T}\}}|y|^{q-2d-\varepsilon}   \,{\mathrm d} y\\
 &=& \bar c^q c_d   {\mathrm e}^{(\Lambda +\frac{1}{2}\sigma^2 q) q T}{\mathrm e}^{-\gamma d T} 
\int_0^{\sqrt{d}{\mathrm e}^{-\gamma T}} r^{q-d-1-\varepsilon}  \,{\mathrm d}  r\\
&=& \bar c^q c_{d,q,\varepsilon} {\mathrm e}^{(\Lambda-\gamma +\frac{1}{2}\sigma^2 q) q T}  {\mathrm e}^{\gamma \varepsilon T}.
\end{eqnarray*}
Therefore, the assumptions of Lemma \ref{GRR} are satisfied for almost all $\omega \in \Omega$ and we obtain
\begin{eqnarray*}
&&{\mathbf P}\{\sup_{x,y \in {\mathcal X}_T} \sup_{0 \le t \le T} \rho(\phi_t(x),\phi_t(y))\ge u\}\\ 
&\le& {\mathbf P}\Bigg\{ V^{1/q}(\omega)  \int_0^{4\sqrt{d}{\mathrm e}^{-\gamma T}} s^{\frac{\varepsilon}{q}-1}  \,{\mathrm d}  s \ge c_{d,q,\varepsilon} 
u  \Bigg\}\\
&\le& {\mathbf E} V {\mathrm e}^{-\gamma \varepsilon T} c_{d,q,\epsilon} u^{-q}\\
&\le& \bar c^q c_{d,q,\varepsilon} {\mathrm e}^{(\Lambda-\gamma+\frac{1}{2}\sigma^2 q)qT} u^{-q}. 
\end{eqnarray*}
Taking logarithms, dividing by $T$, letting $T \to \infty$ and then $\varepsilon \to 0$ and 
optimizing over $q>d$ yields the claim in Theorem \ref{SB}.
\hfill $\Box$\\

Occasionally, the GRR-Lemma is formulated only for $p$ being the identity (e.g.~in \cite{ai}). The following proof shows that 
we don't loose anything in this case but a few modifications are necessary.\\

\noindent{\bf Proof of Theorem \ref{SB} using Lemma \ref{GRR} with $p$=id.}
Fix $u>0$. We start as in the previous proof except that we choose $p(s)=s$, $q>2d$, $Q \in (0,1)$, $qQ \ge 1$, and
\begin{eqnarray*}
V(x,y)&:=&\sup_{0 \le t \le T} \left(\frac{\rho(\phi_t(x),\phi_t(y))\wedge u}{|x-y|}  \right)^q\\
V&:=&\int_{{\mathcal X}_T} \int_{{\mathcal X}_T} V(x,y)  \,{\mathrm d}  x  \,{\mathrm d}  y.
\end{eqnarray*}
Using Chebychev's inequality and (H), we get
\begin{eqnarray*}
{\mathbf E} V(x,y)&=&\int_0^{(u/|x-y|)^q} {\mathbf P}\{ V(x,y) \ge s\}  \,{\mathrm d}  s\\
&\le& {\mathbf E} \sup_{0 \le t \le T} \left(\frac{\rho(\phi_t(x),\phi_t(y))}{|x-y|}  \right)^{qQ} \int_0^{(u/|x-y|)^q} s^{-Q}  \,{\mathrm d}  s\\
&\le&\bar c^{qQ} (1-Q)^{-1} \exp\{(\Lambda + \frac{1}{2}\sigma^2qQ)qQT\} \left(\frac{u}{|x-y|}\right)^{q(1-Q)}.
\end{eqnarray*}
Hence  Lemma \ref{GRR} with $p$=id implies
\begin{eqnarray}
&&{\mathbf P}\{\sup_{x,y \in {\mathcal X}_T} \sup_{0 \le t \le T} \rho(\phi_t(x),\phi_t(y))\ge u\}\label{first}\\ 
&\le&{\mathbf P}\bigg\{V^{1/q} \ge c_{q,d} u {\mathrm e}^{\gamma \left(-\frac{2d}{q} + 1\right)T} \bigg\}\nonumber\\
&\le&{\mathbf E} V c_{d,q} u^{-q} {\mathrm e}^{(2d-q)\gamma T} \nonumber\\
&\le& c_{q,Q,d} \bar c^{qQ}  {\mathrm e}^{(\Lambda +\frac{1}{2}\sigma^2qQ)qQT} {\mathrm e}^{(2d-q)\gamma T} u^{-qQ}
\int_{{\mathcal X}_T} \int_{{\mathcal X}_T} |x-y|^{-q(1-Q)}  \,{\mathrm d}  x   \,{\mathrm d} y.\nonumber
\end{eqnarray}
The double integral is finite if $q(1-Q)<d$. Observe that for any $\kappa>d$ we can find $q>2d$ and $Q \in (0,1)$ such that
$qQ=\kappa$ and  $q(1-Q)<d$. Therefore we obtain the same asymptotics for  (\ref{first}) as in the previous proof. \hfill$\Box$\\

Theorem \ref{SB} can be improved in case $\phi_t$ is a homeomorphism on ${\mathbf R}^d$ for each $t \ge 0$ 
and each $\omega \in \Omega$.

\begin{Corollary}
Let $\phi_t$ be a homeomorphism on ${\mathbf R}^d$, $d \ge 1$ for each $t \ge 0$ and $\omega \in \Omega$. If $\phi$ satisfies {\rm{(H)}} with 
respect to the Euclidean norm $\rho$, then the conclusion of Theorem \ref{SB} holds when in the definition of $I$, $d$ is 
replaced by $d-1$.
\end{Corollary}

\noindent {\bf Proof.} Due to the homeomorphic property, the sup over ${\mathcal X}_T$ in Theorem \ref{SB} is attained on one of the 
faces of ${\mathcal X}$. Applying Theorem \ref{SB} to each of the faces (which have dimension $d-1$), the assertion in the corollary follows.
\hfill $\Box$

\section{Examples and Complements}\label{Examples}

Let us first show that a solution flow of a stochastic differential equation on ${\mathbf R}^d$ with Lipschitz coefficients satisfies 
hypothesis (H') and therefore also (H).

For each $x \in {\mathbf R}^d$, let $t \mapsto M(t,x)$ be an ${\mathbf R}^d$-valued continuous martingale with $M(0,x)=0$ such that the
joint quadratic variation can be represented as
$$
\langle M(.,x.\omega),\, M(.,y,\omega) \rangle_t=\int_0^t a(s,x,y,\omega)  \,{\mathrm d}  s,
$$ 
for a jointly measurable matrix--valued function $a$ which is continuous in $(x,y)$ and predictable in $(s,\omega)$.
Defining
$$
{\cal{A}}(s,x,y,\omega):=a(s,x,x,\omega) - a(s,y,x,\omega) - a(s,x,y,\omega) + a(s,y,y,\omega),
$$
we will require that $a$ satisfies the following Lipschitz property: there exists some constant $a \ge 0$ such that for all 
$x,y \in {\mathbf R}^d$, all $s \ge 0$ and almost all $\omega$, we have
$$
\|{\cal{A}}(s,x,y,\omega)\| \le a^2 |x-y|^2,
$$
where $\|.\|$ denotes the operator norm. Note that
$$
{\cal{A}}(t,x,y,\omega)=\frac{{\mathrm d} }{{\mathrm d}  t} \langle M(.,x) - M(.,y) \rangle_t.
$$
Further, we assume that $b: [0,\infty) \times {\mathbf R}^d \times \Omega \to {\mathbf R}^d$ is a vector field which is jointly measurable, 
predictable in $(t,\omega)$ and Lipschitz continuous with constant $b$ in the spatial variable uniformly in $(t,\omega)$. 
In addition, we require, that the functions $a(.)$ and $b(.)$ are bounded on each compact subset of 
$[0,\infty) \times {\mathbf R}^{d \times d}$ resp. $[0,\infty) \times 
{\mathbf R}^d$ uniformly w.r.t.~$\omega \in \Omega$. 
Under these assumptions, it is well known that the {\em Kunita type} stochastic differential equation
\begin{equation}\label{stochdiff}
{\mathrm d} X(t)=b(t,X(t))  \,{\mathrm d}  t + M({\mathrm d}  t, X(t))
\end{equation}
generates a stochastic flow of homeomorphisms $\phi$ (see \cite{k}, Theorem 4.5.1), i.e.
\begin{itemize}
\item[i)] $t \mapsto \phi_{s,t}(x), \, t \ge s$ solves (\ref{stochdiff}) with initial condition $X(s)=x$ for all 
$x \in {\mathbf R}^d$, $s\ge 0$. 
\item[ii)] $\phi_{s,t}(\omega)$ is a homeomorphism on ${\mathbf R}^d$ for all $0 \le s \le t$ and all $\omega \in \Omega$.
\item[iii)] $\phi_{s,u}=\phi_{t,u} \circ \phi_{s,t}$ for all $0 \le s \le t \le u$  and all $\omega \in \Omega$.
\item[iv)] $(s,t,x) \mapsto \phi_{s,t}(x)$ is continuous. 
\end{itemize}
 We will write $\phi_t(x)$ instead of $\phi_{0,t}(x)$.

For readers who are unfamiliar with Kunita type stochastic differential equations, we point out that if one replaces 
the term $M({\mathrm d} t,X(t))$ in equation (\ref{stochdiff}) by $\sum_{i=1}^m \sigma_i (X(t)) \, {\mathrm d} W_i(t)$, 
where $W_i$ are independent scalar standard Brownian motions and the functions $\sigma_i: {\mathbf R}^d \to 
 {\mathbf R}^d$ are Lipschitz continuous, then the Lipschitz condition imposed above holds. In fact
$$
{\cal{A}}(t,x,y,\omega)=\sum_{i=1}^m (\sigma_i(x)-\sigma_i(y))(\sigma_i(x)-\sigma_i(y))^T.
$$

The following Lemma is identical with Lemma 5.1 in \cite{css2}. The proof below is slightly more elementary since it 
avoids the use of a comparison theorem by Ikeda and Watanabe.
\begin{Lemma}\label{sde}
Under the assumptions above, \textnormal{(H')} holds with $\sigma=a$ and $\Lambda=b+(d-1)a^2/2$.
\end{Lemma}

\noindent {\bf Proof.} Fix $x,\,y \in {\mathbf R}^d$, $x \neq y$ and define

$$D_{t} := \phi_t(x) - \phi_t(y),\;\;Z_{t} := \frac{1}{2} \log(|D_t|^2). $$
Therefore, $Z_t = f (D_t)$ where $f(z) := \frac{1}{2}\log(|z|^2)$. Note that $D_t \neq 0$ for all $t\ge0$ by the homeomorphic
property. Using It\^{o}'s formula, we get
\begin{eqnarray*}
   {\mathrm d}  Z_t &=& \frac{D_{t} \cdot \left ( M({\mathrm d}  t, 
  \phi_t (x)) - M({\mathrm d} t,\phi_t (y)) \right )}{| D_{t}|^2 } + 
  \frac{D_t \cdot \left ( b (t, \phi_t(x)) - b(t, 
  \phi_t (y)) \right )}{| D_{t}|^2}  \,{\mathrm d}  t\\
 & & + \frac{1}{2} \frac{1}{|D_{t}|^2} {\mbox {Tr}} \left ( 
  {\cal A} (t, \phi_t(x), \phi_t(y), \omega) \right )  {\mathrm d}  t\\
 & & - \sum_{i,j} \frac{D^i _{t} D^j _{t}}{(|D_{t}|^2)^2} 
  {\cal A}_{i,j} (t, \phi_t(x), \phi_t(y), \omega)  \,{\mathrm d}  t.
\end{eqnarray*}
We define the local martingale $N_{t}, t \ge 
0$ by
$$
N_t = \int^t_{0}\frac{D_{s}}{|D_{s}|^2} \cdot \left ( M({\mathrm d}  s, \phi_s(x)) - M({\mathrm d}  s, 
\phi_s(y))\right )
$$
and obtain
$$
Z_t = Z_{0} + N_t + \int^t_{0} \alpha (s,\omega)  \,{\mathrm d}  s,  
$$
where
$$
\sup_{x,y}\sup_s {\mathrm {esssup}}_{\omega}  
|\alpha (s,\omega)|  \le b + (d-1)a^2/2 =:\Lambda
$$ 
and
\begin{equation}\label{kappa}
{\mathrm d}  \langle N \rangle_{t} = \sum_{i,j} 
\frac{D^i_{t}D^j_{t}}{(|D_{t}|^2)^2}{\cal{A}}_{i,j}(t, 
\phi_t (x), \phi_t (y), \omega)  \,{\mathrm d}  t \le a^2   \,{\mathrm d} t. 
\end{equation}
Since $N$ is a continuous local martingale with 
$N_{0} = 0$, there exists a standard Brownian motion $W$ (possibly on an enlarged probability space) 
such that $N_t=a W_{\tau(t)}$, $t \ge 0$ and (\ref{kappa}) 
implies $\tau(t) \le t$ for all $t \ge 0$. Hence
\begin{equation}
Z_{t} \le \log |x-y| + a \; W^{*}_{t} + \Lambda t. 
\end{equation}
Exponentiating the last inequality completes the proof of the lemma. \hfill $\Box$\\

The following simple example shows that the upper bound in Theorem \ref{SB} is sharp for $\gamma \ge \Lambda + \sigma^2 d$.\\

\noindent{\bf Example}
Consider the linear stochastic differential equation
$$
\,{\mathrm d}   X(t) = (\Lambda  + \frac{1}{2} \sigma^2)\,X(t)\, {\mathrm d} t + \sigma X(t)\, {\mathrm d} W(t), \;X(0)=x \in {\mathbf R}^d,
$$
where $W(t),\,t \ge 0$ is a one-dimensional Brownian motion, $\Lambda \ge 0$ and $\sigma > 0$. 
The solution (flow) $\phi_t(x)$ is given by
$$
\phi_t(x) = x {\mathrm e}^{\Lambda t + \sigma W(t)},
$$
which satisfies (H') and hence (H). If ${\mathcal X}$ is a cube of side length ${\mathrm e}^{-\gamma T}$ in ${\mathbf R}^d$ for some $\gamma \ge \Lambda$
and $u>0$, then 
$$
\sup_{x,y \in {\mathcal X}} |\phi_t(x)-\phi_t(y)|=\sqrt{d}{\mathrm e}^{-\gamma T} {\mathrm e}^{\Lambda t + \sigma W(t)} 
$$
and
$$
\limsup_{T \to \infty} \frac{1}{T} \log {\mathbf P}\{\sup_{x,y \in {\mathcal X}} \sup_{0 \le t \le T} |\phi_t(x)-\phi_t(y)| \ge u\} 
= -\frac{1}{2 \sigma^2}(\gamma-\Lambda)^2,
$$
for $\gamma \ge \Lambda$ and $u>0$. \hfill $\Box$\\

Next, we provide an example which shows that the conclusion in Theorem \ref{SB} is sharp also for 
$\gamma < \Lambda + \sigma^2 d$.\\

\noindent{\bf Example} Let $h:{\mathbf R}^d \to [0,\infty)$ be Lipschitz continuous with Lipschitz constant 2 and support contained in 
$[-1/2,1/2]^d$. Further suppose that $h(0)=1$. Let $W^i$, $i \in {\mathbf Z}^d$ be independent standard Brownian motions and let $\Lambda \ge 0$ 
and $\sigma>0$ be constants. For $\delta>0$, define
$$
\phi_t (x) := \sum_{i \in {\mathbf Z}^d} \delta h\left( \frac{x}{\delta}-i \right) 
{\mathrm e}^{\Lambda t + \sigma W^i_t},\;x \in {\mathbf R}^d,\,t \ge 0.
$$
Note that at most one term in the sum is nonzero. Therefore
$$
\left( {\mathbf E} \sup_{0 \le t \le T} |\phi_t(x)-\phi_t(y)|^q \right)^{1/q} \le 2|x-y|  {\mathrm e}^{(\Lambda + \frac{1}{2}q\sigma^2)T}
$$
for each $q \ge 1$, so (H) is satisfied with $\bar c=2$.
Let $T>0$, $\gamma > \Lambda$, ${\mathcal X}=[0,{\mathrm e}^{-\gamma T}]^d$ and $\delta={\mathrm e}^{-\xi T}$, where $\xi>\gamma$ will be optimized later.
Since the processes $\phi_t(i\delta)$, $i \in {\mathbf Z}^d$ are independent and identically distributed, we conclude
\begin{align*}
&{\mathbf P}\{\sup_{x \in {\mathcal X}} \phi_T(x) \le \delta {\mathrm e}^{\Lambda T} + 1\}\\ 
&\le {\mathbf P}\{\max_{i \in {\mathbf Z}^d, i\delta \in {\mathcal X}} \phi_T(i\delta) \le \delta {\mathrm e}^{\Lambda T} + 1\}\\
&\le \left({\mathbf P}\{\phi_T(0) \le \delta {\mathrm e}^{\Lambda T} + 1\}\right)^{\exp\{(\xi-\gamma)dT\}}\\
&= \left( 1-{\mathbf P}\left\{W_1 > \frac{1}{\sigma \sqrt{T}} \log (1+{\mathrm e}^{(\xi-\Lambda)T})\right\}\right)^{\exp\{(\xi-\gamma)dT\}}.
\end{align*}
 From this and the asymptotic behaviour of the last probability, it follows that the last term will converge to 0 as $T \to \infty$ 
provided that
$$
2\sigma^2 d(\xi-\gamma)>(\xi-\Lambda)^2,
$$
which holds true in case $\xi=\Lambda + \sigma^2 d$ and $\gamma \in (\Lambda,\Lambda + \sigma^2 d/2)$. Since the probability that 
the infimum of $\phi_T(x),\,x \in {\mathcal X}$ is at most $\delta \exp\{\Lambda T\}$ converges to one as $T \to \infty$, we obtain
for $\gamma < \Lambda + \sigma^2 d/2$
$$
\limsup_{T \to \infty} \frac{1}{T} \log {\mathbf P}\{\sup_{x,y \in {\mathcal X}} \sup_{0 \le t \le T} |\phi_t(x)-\phi_t(y)| \ge 1\} = 0
$$
(in fact we just showed that this is true even if the sup over $t$ is replaced by $T$). 
Similarly, we obtain 
\begin{equation}\label{uups}
\limsup_{T \to \infty} \frac{1}{T} \log {\mathbf P}\{\sup_{x \in {\mathcal X}} \phi_T(x) \ge \delta {\mathrm e}^{\Lambda T} + 1\} 
= (\xi-\gamma)d-\frac{1}{2 \sigma^2}(\xi-\Lambda)^2
\end{equation}
in case the last expression is strictly negative which holds true in case 
$\gamma \in (\Lambda + \frac{1}{2}\sigma^2 d,\Lambda + \sigma^2 d)$ and $\xi = \Lambda +  \sigma^2 d$. Inserting this value for 
$\xi$ in (\ref{uups}) yields
$$
\limsup_{T \to \infty} \frac{1}{T} \log {\mathbf P}\{\sup_{x \in {\mathcal X}} \phi_T(x) \ge \delta {\mathrm e}^{\Lambda T} + 1\} =-I(\gamma)
$$
for all $\gamma < \Lambda+\sigma^2 d$ with $I(\gamma)$ defined as in Theorem \ref{SB}.\\

The reader may complain, that in this example the field $\phi$ actually depends on $T$ (via $\delta=\delta(T)$),
i.e.~as we let $T \to \infty$, we keep changing $\phi$. It is easy to see however, that we can define 
a single field $\phi$ by spatially  piecing together fields as above for an appropriate sequence $T_i \to \infty$.\\  

\noindent {\bf Remark} \hspace{.3cm} The previous example(s) show that the conclusion in Theorem \ref{SB} 
is sharp, but in the last example $\phi$ is not a stochastic flow of homeomorphisms. Can we do better in 
that case? The following theorem shows that we can, provided the flow is $C^1$. 
More precisely, we consider a stochastic flow of homeomorphisms $\phi$ as introduced at the beginning of this 
section and require that it has -- in addition -- independent and stationary increments and that its law 
is invariant under shifts in ${\mathbf R}^d$. We will call such a flow a {\em translation invariant Brownian flow}.

\begin{Theorem}\label{diff} 
Let $\phi$ be a translation invariant Brownian flow on ${\mathbf R}^d$ such that the map 
$(t,x) \mapsto D\phi_t(x)$ is continuous (for all $\omega \in \Omega$). In addition, we 
assume that there exist $\bar c \ge 1, \Lambda \ge 0$ and $\sigma > 0$ 
and a standard Wiener process such that for each $T \ge 0$, we have
\begin{equation}\label{Voraussetzung}
\|D\phi_T(0)\| \le \bar c \exp\{\sup_{0 \le s \le T} (\Lambda s + \sigma W_s)\}.
\end{equation}
Then, for each $u>0$ and $\xi \ge 0$, we have
$$
\limsup_{T \to \infty} \frac{1}{T} \log {\mathbf P} \left\{ \sup_{|x| \le \exp\{-(\Lambda+\xi)T\}} \sup_{0 \le t \le T} 
|\phi_t(x)-\phi_t(0)| \ge u \right\} \le -\frac{\xi^2}{2 \sigma^2}.\\
$$
\end{Theorem}

Note that, due to the fact that the flow is translation invariant and stationary, the statement is invariant 
under a shift in space and time as well.
We mention that the hypotheses of the theorem are for example fulfilled for isotropic Brownian flows, see \cite{bh}.\\

\noindent {\bf Proof of Theorem \ref{diff}.}
Fix $\xi>0$, $\varepsilon \in (0,1)$, $z>\varepsilon + \bar c \exp\{-\frac{\sigma^2}{\xi}\log (1-\varepsilon)\}$ 
and $u>0$. We abbreviate $D_t:=\|D\phi_t(0)\|$. Let
\begin{eqnarray*}
\tau_z&:=&\inf \{t \ge 0: \Lambda t + \sigma W_t \ge \log \frac{z}{\bar c}\}.
\end{eqnarray*}
Using the formula for the Laplace transform of the hitting time of Brownian motion with drift (\cite{bs}, 
page 223, formula 2.2.0.1), we get for $\lambda > 0$
\begin{eqnarray}\label{absch}
{\mathbf E} {\mathrm e}^{-\lambda \tau_z} 
&=&\exp\{ \frac{1}{\sigma^2}  (\Lambda - \sqrt{2 \lambda \sigma^2 + \Lambda^2}) \log (z/\bar c)\}.
\end{eqnarray}
For $\delta > 0$, we define
$$
\tilde T^{(\delta)}:=\inf\{t>0: \sup_{|x| \le \delta} |\phi_t(x)-\phi_t(0)| \ge \delta z\}.
$$
Since the flow $\phi$ is $C^1$ and $\tau_{z-\varepsilon} < \infty$, there exists $\delta_0=\delta_0 (z,\varepsilon) >0$ such that
$$
{\mathbf P}\bigg\{ \sup_{0 \le t \le \tau_{z-\varepsilon}}\sup_{|x|\le \delta} \frac{|\phi_t(x)-\phi_t(0)|}{\delta} 
< \sup_{0 \le t \le \tau_{z-\varepsilon}} D_t + \varepsilon\bigg\} \ge 1-\varepsilon
$$
for all $\delta \in (0,\delta_0]$. Note that (\ref{Voraussetzung}) implies 
$\sup_{0 \le t \le \tau_{z-\varepsilon}} D_t+\varepsilon \le z$. Hence
\begin{equation}\label{laplace}
{\mathbf E}{\mathrm e}^{-\lambda \tilde T^{(\delta)}} \le {\mathbf E}{\mathrm e}^{-\lambda \tau_{z-\varepsilon}} + \varepsilon
\end{equation}
for  $\delta \in (0,\delta_0]$ and all $\lambda > 0$. Define $\hat u:= u \wedge \delta_0$.

Let $T>0$ such that $\exp\{-(\Lambda + \xi)T\} < \hat u$. Further,
let $T_1,T_2,...$ be independent random variables such that the laws of $T_j$ and $\tilde T^{(\delta_j)}$ coincide, 
where $\delta_j= \exp\{-(\Lambda + \xi )T\} z^{j-1}$. Define
\begin{equation}\label{ceil}
m= \left\lfloor \frac{(\Lambda+\xi)T + \log \hat u}{\log z} \right\rfloor.
\end{equation}
Using the fact that $\phi$ has independent and stationary increments and Markov's inequality, we obtain
\begin{align*}
{\mathbf P}\{ \sup_{|x| \le \exp\{-(\Lambda+\xi)T\}} \sup_{0 \le t \le T} 
&|\phi_t(x)-\phi_t(0)|    \ge u\} \le {\mathbf P} \{ \sum_{j=1}^m T_j \le T\}\\
&= {\mathbf P}\{ \exp\{ -\lambda  \sum_{j=1}^m T_j\} \ge \exp\{- \lambda T\}\}\\
&\le \exp\{\lambda T\} \max_{j=1,...m} ({\mathbf E} \exp\{ -\lambda T_j\})^m\\
&\le \exp\{\lambda T\} ({\mathbf E} \exp\{ -\lambda \tau_{z-\varepsilon}\} + \varepsilon)^m,  
\end{align*}
where we used (\ref{laplace}) and $\hat u \le \delta_0$ in the last step. 
Using  (\ref{absch}) and (\ref{ceil})  
and inserting  $\lambda:=\frac{1}{2 \sigma^2} ((\Lambda+\xi)^2-\Lambda^2)$, we get
\begin{align*}
&\log {\mathbf P}\{ \sup_{|x| \le \exp\{-(\Lambda+\xi)T\}} \sup_{0 \le t \le T} 
|\phi_t(x)-\phi_t(0)|    \ge u\}\\ 
&\hspace{.3cm} \le \frac{(\Lambda + \xi)^2-\Lambda^2}{2 \sigma^2} T + \left\lfloor \frac{(\Lambda+\xi)T + \log \hat u}{\log z}
\right\rfloor \log\left( \exp\{-\frac{\xi}{\sigma^2}\log \frac{z-\varepsilon}{\bar c}\} + \varepsilon \right).   
\end{align*}
Dividing by $T$, and letting (in this order) $T \to \infty$, $\varepsilon \to 0$ and $z \to \infty$, 
the assertion follows. \hfill $\Box$

\section{Dispersion of Sets: Upper Bounds}

We will now formulate the dispersion result mentioned in the introduction and prove it using Theorem \ref{SB}. In addition to 
hypothesis (H) we require a growth condition for the one-point motion. In Proposition \ref{1pt} we will provide 
explicit conditions on the coefficients of a stochastic differential equation which guarantee that the associated stochastic 
flow fulfills that condition. The value of the linear bound $K$ in Theorem \ref{diam} improves previous ones in 
\cite{css2,ls01,ls03} but the main improvement is its simpler proof.

\begin{Theorem}\label{diam}
Let $\phi: [0,\infty)\times {\mathbf R}^d \times \Omega \to {\mathbf R}^d$ be a continuous random field satisfying
\begin{itemize}
\item[\textnormal{(i)}] \textnormal{(H)}.
\item[\textnormal{(ii)}] There exist $A>0$ and $B \ge 0$ such that for each $k>0$ and each bounded set 
$S \subset {\mathbf R}^d$, we have
$$
\limsup_{T \to \infty} \frac{1}{T} \log \sup_{x \in S}  {\mathbf P}\left\{ \sup_{0 \le t \le T} |\phi_t(x)| \ge kT\right\} \le 
-\frac{(k-B)_+^2}{2A^2},
$$
where $r_+=r \vee 0$ denotes the positive part of $r \in {\mathbf R}$.
\end{itemize}
Let ${\mathcal X}$ be a compact subset of ${\mathbf R}^d$ with box (or upper entropy) dimension $\Delta > 0$. Then
\begin{equation}\label{LB}
\limsup_{T \to \infty} \left( \sup_{t \in [0,T]} \sup_{x \in {\mathcal X}} \frac{1}{T}| \phi_t(x)| \right) \le K \:a.\,s.,
\end{equation}
where
\begin{eqnarray}\label{gen}
K&=& \left\{
\begin{array}{ll}
B+A\sqrt{ 2\Delta  \left(\Lambda + \sigma^2\Delta +\sqrt{ \sigma^4\Delta^2 + 2\Delta \Lambda \sigma^2} \right) 
} & \mbox{ if }\; \Lambda \ge \Lambda_0\\
B+A\sqrt{2\Delta \frac{d}{d-\Delta}\left( \Lambda+\frac{1}{2}\sigma^2d\right)}& \mbox{ otherwise}\;,
\end{array}
\right.
\end{eqnarray}
where 
$$
\Lambda_0:=\frac{\sigma^2 d}{\Delta} \left(\frac{d}{2} -\Delta \right).
$$
\end{Theorem}

\noindent{\bf Proof.} 
Let $N({\mathcal X},r)$, $r>0$ denote the minimal number of subsets of ${\mathbf R}^d$ of diameter at most $r$ 
which cover ${\mathcal X}$. By definition, we have
$$
\Delta=\limsup_{r \downarrow 0} \frac{\log N({\mathcal X},r)}{\log \frac{1}{r}}.
$$
Choose $\varepsilon>0$ and $r_0>0$ such that $\log N({\mathcal X},r) \le (\Delta+\varepsilon)\log \frac{1}{r}$ for all $0 <r\le r_0$. 
Further, let $\gamma,T>0$ satisfy ${\mathrm e}^{-\gamma T}\le r_0$. Then $N({\mathcal X},{\mathrm e}^{-\gamma T}) \le 
\exp\{\gamma T(\Delta +\varepsilon)\}$. Let 
${\mathcal X}_i,\, i=1,\dots,N({\mathcal X},{\mathrm e}^{-\gamma T})$ be compact sets of diameter at most ${\mathrm e}^{-\gamma T}$ which cover 
${\mathcal X}$ and 
choose arbitrary points $x_i \in {\mathcal X}_i$. Define 
$$\widetilde{{\mathcal X}}:=\{x_i,\,i=1,\dots,N({\mathcal X},{\mathrm e}^{-\gamma T})\}.$$ For 
$\kappa>0$, we have 
$$
{\mathbf P}\{\sup_{x \in {\mathcal X}} \sup_{0 \le t \le T} | \phi_t(x)| \ge \kappa T \} \le S_1 + S_2,
$$
where
$$
S_1:=\exp\{\gamma T(\Delta +\varepsilon) \} \max_{x \in \widetilde{{\mathcal X}}} 
{\mathbf P}\{\sup_{0 \le t \le T} | \phi_t(x)| \ge \kappa T -1\}
$$
and
$$
S_2:=\exp\{\gamma T(\Delta+\varepsilon) \} \max_{i} {\mathbf P}\{\sup_{0 \le t \le T} {\mathrm{diam}} (\phi_t({\mathcal X}_i)) \ge 1\}.
$$
Using (ii) in the theorem, we get
\begin{equation}\label{S1}
\limsup_{T \to \infty} \frac{1}{T} \log S_1 \le \gamma(\Delta+\varepsilon)-\frac{(\kappa-B)_+^2}{2A^2}.
\end{equation}
Further, Theorem \ref{SB} implies
\begin{equation}\label{S2}
\limsup_{T \to \infty} \frac{1}{T} \log S_2 \le \gamma(\Delta+\varepsilon)-I(\gamma).
\end{equation}
Therefore,
\begin{eqnarray*}
\zeta(\gamma,\kappa)&:=&\limsup_{T \to \infty} \frac{1}{T} \log {\mathbf P}\{\sup_{x \in {\mathcal X}} 
\sup_{0 \le t \le T} | \phi_t(x)| \ge \kappa T\}\\ 
&\le&\gamma \Delta - \left( \frac{(\kappa-B)_+^2}{2A^2} \wedge I(\gamma) \right).
\end{eqnarray*}
Let $\gamma_0$ be the unique positive solution of $I(\gamma)=\gamma \Delta$, where $I(\gamma)$ is defined in Theorem \ref{SB}. Then 
\begin{eqnarray*}
\gamma_0&=& \left\{
\begin{array}{ll}
\frac{d}{d-\Delta} (\Lambda + \frac{1}{2} \sigma^2d)   
& \mbox{ if }\; \Lambda \le \frac{\sigma^2 d}{\Delta} \left( \frac{d}{2} -\Delta \right)\\
\Lambda+\sigma^2\Delta+\sqrt{2\Lambda\sigma^2\Delta+\sigma^4\Delta^2} & \mbox{ otherwise}
\end{array}
\right.
\end{eqnarray*}
and $\zeta(\gamma,\kappa)<0$ whenever $\gamma>\gamma_0$ and $\kappa>\kappa_0(\gamma)$, where
$$
\kappa_0(\gamma):= B+A\sqrt{2\gamma \Delta}.
$$
Therefore, for any $\gamma>\gamma_0$ we have 
$$
\sum_{n=1}^{\infty} {\mathbf P} \left\{\frac{1}{n} \sup_{x \in {\mathcal X}} \sup_{0 \le t \le n} | \phi_t(x) | \ge 
B+A\sqrt{2\gamma \Delta}\right\} < \infty.
$$
Using the Borel--Cantelli  Lemma, we obtain
$$
\limsup_{T \to \infty} \frac{1}{T} \sup_{x \in {\mathcal X}} \sup_{0 \le t \le T} | \phi_t(x) | \le 
K:=B+A\sqrt{2\gamma_0 \Delta} \; 
\mbox{ a.s.}
$$
which proves the theorem. \hfill $\Box$\\

\begin{Corollary}
Assume in addition to the hypotheses in Theorem \ref{diam} that $\phi_t$ is a (random) homeomorphism on ${\mathbf R}^d$ for each $t \ge 0$. 
Then (\ref{LB}) holds with $\Delta$ replaced by $\Delta \wedge(d-1)$.
\end{Corollary}

\noindent {\bf Proof.} Let ${\mathcal X}$ be compact and have box dimension $>(d-1)$ and let $\tilde {\mathcal X}$ be a compact set which contains 
${\mathcal X}$ such that $\partial \tilde {\mathcal X}$ has box dimension $d-1$. We then apply Theorem \ref{diam} to $\partial \tilde {\mathcal X}$ instead of 
${\mathcal X}$. By the homeomorphic property, we know that $\sup_{x \in {\mathcal X}} |\phi_t(x)| \le\sup_{x \in \partial \tilde {\mathcal X}} |\phi_t(x)|$ 
and the assertion of the corollary follows. \hfill $\Box$ \\

\noindent {\bf Remark} \hspace{.3cm} If $\phi$ is a flow which satisfies the assumptions of Theorem \ref{diff}, then the upper 
bound for $K$ in Theorem \ref{diam} can be improved by changing $I(\gamma)$ in (\ref{S2}) accordingly. In this case 
the upper formula for $K$ in (\ref{gen}) holds for all values of $\Lambda$.\\

Now we provide a class of stochastic differential equations for which the assumptions of the previous theorem are satisfied.
For simplicity we will assume that the drift $b$ is autonomous and deterministic (if it is not, but the bound on $b$ 
in the proposition is uniform with respect to $(t,\omega)$, then the proposition and its proof remain true without further 
change).

\begin{Proposition}\label{1pt}
Let the assumptions of Lemma \ref{sde} be satisfied and assume in addition (for simplicity) that $a(.)$ and $b(.)$ 
are deterministic and autonomous. Further, we require that there exists $A>0$ such that 
$\|a(x,x)\| \le A^2$ for all $x$ and that 
$$
\limsup_{|x| \to \infty}\, \frac{x}{|x|} \cdot b(x) \le B \in {\mathbf R}.
$$
Then for each compact set $S$ and each $k > 0$
\begin{eqnarray*}
\limsup_{T \to \infty} \frac{1}{T} \log \sup_{x \in S} {\mathbf P} \left\{ \sup_{0 \le t \le T} |\phi_t(x)| \ge kT\right\} \le \left\{
\begin{array}{ll}
-\frac{(k-B)_+^2}{2A^2}
& \mbox{ if } k \ge -B\\
2Bk \frac{1}{A^2} & \mbox{ otherwise.}
\end{array}
\right.\\
\end{eqnarray*}
\end{Proposition}
\noindent {\bf Proof.} 
Let $S$ be a compact subset of ${\mathbf R}^d$ and $k > B $ (otherwise there is nothing to show).
Fix $0<\varepsilon<k-B$ and let $r_0> 1$ be such that 
$$
\frac{x}{|x|} \cdot b(x) + \frac{d-1}{2|x|} A^2 \le B+\varepsilon \;\;\mbox{ for all } |x| \ge r_0
$$
and such that $S$ is contained in a ball around 0 of radius $r_0$.
Let $h$ be an even smooth function from ${\mathbf R}$ to ${\mathbf R}$ such that 
$h(y)=|y|$ for $|y| \ge 1$ and $|h'(y)| \le 1$ for all 
$y \in {\mathbf R}$ and  
define $\rho_t(x)=h(|\phi_t(x)|)$. Applying It\^o's formula, we get
\begin{align*}
{\mathrm d} \rho_t(x)={\mathrm d} N_t + f(\phi_t(x))\,{\mathrm d} t,
\end{align*}
where 
\begin{eqnarray*}
N_t&=& \sum_{i=1}^d \int_0^t h'(\rho_s(x))\frac{\phi^i_s(x)}{\rho_s(x)} M^i({\mathrm d} s,\phi_s(x)) \,\,\,\,\,\,\,\,\text{ and} \\
f(x)&=& \frac{x}{|x|}\cdot b(x) + \frac{1}{2|x|}{\mathrm{Tr}}\,a(x,x)-\frac{1}{2|x|^3}x^Ta(x,x)x\\
&\le& \frac{x}{|x|}\cdot b(x) + \frac{d-1}{2|x|} A^2
\le B + \varepsilon \;\;\text{ on } \{|x|\ge  r_0\}.
\end{eqnarray*}
For the quadratic variation of $N$, we have the following bound:
\begin{align*}
\langle N \rangle_t - \langle N \rangle_s \le \int\limits_s^t \frac{1}{\rho_u^2(x)}\phi_u^T(x)
a(\phi_u(x),\phi_u(x))\phi_u(x)\,{\mathrm d} u\le A^2 (t-s)\,.
\end{align*}
The continuous local martingale  $N$  can be represented (possibly on an 
enriched probability space) in the form  $N_t=A\, W_{\tau(t)}$, 
where $W$ is a standard Brownian motion and the family of stopping times $\tau(s)$ satisfies 
$\tau(t)-\tau(s) \le t- s $ whenever $s \le t$.
For $|x| \le r_0 < kT$ we get
\begin{align*}
&{\mathbf P} \left\{ \sup_{0 \le t \le T} |\phi_t(x)| \ge kT\right\}\\ 
&\le {\mathbf P} \left\{\exists \,0\le s \le t \le T: \rho_t(x) - \rho_s(x) \ge kT - r_0 
,\, \inf_{s \le u \le t} \rho_u(x) \ge r_0 \right\}\\
&\le {\mathbf P} \left\{\exists \,0\le s \le t \le T: A(W_{\tau(t)}-W_{\tau (s)} ) + (B + \varepsilon)(t-s) 
\ge kT-r_0 \right\}=:\bar {\mathbf P}
\end{align*}
Now we distinguish between two cases:\\

\noindent {\bf Case 1:} $B \ge 0$. Then 
$$
\bar {\mathbf P} \le {\mathbf P}\left\{ \max_{0 \le s \le 1} W_s -\min_{0 \le s \le 1} W_s \ge 
\frac{k-B-\varepsilon}{A} \sqrt{T} - \frac{r_0}{A\sqrt{T}} \right\}.
$$
The density of the range $R:=\max_{0 \le s \le 1} W_s -\min_{0 \le s \le 1} W_s$ equals
$$
8\sum_{j=1}^{\infty} (-1)^{j-1} j^2 \varphi(jr),
$$
on $[0,\infty)$ (see \cite{f}), where $\varphi$ denotes the density of a standard normal law.
Therefore, for all $u \ge 0$,
$$
{\mathbf P} \{R \ge u \} \le 8 \sum_{j=1}^{\infty}  j \frac{1}{2} \exp\{-\frac{1}{2} j^2 u^2\} 
\sim 4\exp\{-\frac{u^2}{2} \}.
$$
Hence
$$
\limsup_{T \to \infty} \frac{1}{T} \log \sup_{x \in S} {\mathbf P} \left\{ \sup_{0 \le t \le T} |\phi_t(x)| \ge kT\right\} \le 
-\frac{(k-B)_+^2}{2A^2}\,.
$$
\noindent {\bf Case 2:} $B < 0$. We may assume that $\varepsilon >0$ is so small that  also 
$-\tilde B:=(B+\varepsilon)/A<0$.
We have
\begin{equation}\label{LD}
\bar {\mathbf P} \le {\mathbf P} \left\{\exists \,0\le s \le t \le 1: \frac{1}{\sqrt{T}}(W_t-W_s) - \tilde B (t-s) \ge k/A-\frac{r_0}{AT} \right\}.
\end{equation}
To estimate this term, we use large deviations estimates for the standard Wiener process. Let
$$
M:=\{f \in C[0,1]: \, \exists \, 0 \le s \le t \le 1: f_t-f_s - \tilde B(t-s) \ge k/A\}.
$$
The set $M$ is closed in $C[0,1]$ and therefore Schilder's Theorem (\cite{dembo}) implies 
\begin{equation}\label{inf}
\limsup_{T \to \infty} \frac{1}{T} \log {\mathbf P}\{T^{-1/2} \, W \in M \} \le -\inf_{f \in M} I(f),
\end{equation}
where 
\begin{eqnarray*}
I(f)&:=& \left\{
\begin{array}{ll}
\frac{1}{2} \int_0^1 (f_u')^2 \,{\mathrm d} u 
& \mbox{ if } f \mbox{ is absolutely continuous with } L^2 \mbox{ derivative} \\
+\infty & \mbox{ otherwise.}
\end{array}
\right.\\
\end{eqnarray*}
The infimum in (\ref{inf}) can be computed explicitly. 
Let
\begin{eqnarray*}
I&:=& \left\{
\begin{array}{ll}
\frac{1}{2}\left( \frac{k}{A} + \tilde B \right)^2
& \mbox{ if }\; \tilde B \le \frac{k}{A}\\
2\tilde B \frac{k}{A} & \mbox{ otherwise.}
\end{array}
\right.\\
\end{eqnarray*}
For $f_t=(\tilde B+\tilde B \vee (k/A))t$ on $[0,k/((A \tilde B)\vee k)]$ and $f$ constant  
on $[ k/((A \tilde B)\vee k),1]$, we have $f \in M$ and $I(f)=I$.
On the other hand, if $f \in M$ with $I(f)<\infty$, then there exist $0 \le s < t \le 1$ such that 
$f_t-f_s-\tilde B (t-s)\ge k/A$. It follows that 
\begin{align*}
&I(f) \ge \frac{1}{2} \int_s^t (f'_u)^2 \,{\mathrm d} u \ge \frac{1}{2}\frac{1}{t-s}\left(\int_s^t f'_u \,{\mathrm d} u\right)^2\\
&=\frac{1}{2}\frac{1}{t-s} (f_t-f_s)^2 \ge \frac{1}{2}\frac{1}{t-s} \left(\frac{k}{A} + \tilde B (t-s)\right)^2
\ge I.
\end{align*}
Therefore, using (\ref{LD}) and (\ref{inf}), we obtain
$$
\limsup_{T \to \infty} \frac{1}{T} \log \sup_{x \in S} {\mathbf P} \left\{ \sup_{0 \le t \le T} |\phi_t(x)| \ge kT\right\} \le -I
$$
and the proof of the proposition is complete. \hfill $\Box$ \\

\noindent {\bf Remark.} One can modify Theorem \ref{diam} in such a way that it also applies to solution flows generated 
by stochastic differential equations like in the previous proposition with negative $B$. In this case condition (ii) 
in Theorem \ref{diam} has to be changed accordingly. The corresponding linear upper bound will still be strictly positive 
no matter how small $B<0$ is (namely $\gamma_0 \Delta A^2/(-2B)$ as long as this number is at most $-B$). 
In reality however, the linear growth rate turns out to be zero when $B$ is sufficiently small. 
This is shown in \cite{ds}.

\section{Appendix: Proofs of the Chaining Lemmas}
In this section, we provide proofs of those chaining lemmas which are not available in the literature 
in  the form presented here.\\ 

\noindent {\bf Proof of Lemma \ref{Kolmo}.} We skip the proof of the existence of a continuous modification which 
can be found in many textbooks (e.g. \cite{kal}) and only show the estimates, assuming continuity of $Z$.

For $n\in {\bf N}$ define
\begin{eqnarray*}
D_n &:=& \{(k_1, \ldots , k_d)\cdot 2^{-n}; \; k_1, \ldots k_d \in \{1, \ldots , 
2^n  \}  \}\\
\xi_n(\omega) &:=& \max\{\hat \rho(Z_x(\omega), Z_y(\omega)): x, y\in D_n, |x-y|=2^{-n}
\}.
\end{eqnarray*}
The $\xi_n, n\in {\bf N}$ are measurable since $(\hat E, \hat \rho)$  
is separable. Further,
$$
|\{x,y \in D_n : |x-y| = 2^{-n}  \}| \le d\cdot 2^{dn}.
$$
Hence, for  $\kappa \in (0, \frac{b}{a})$,
\begin{eqnarray*}
&&{\bf E} \left( \sum\limits_{n=1}^\infty (2^{\kappa n} \xi_n)^a \right)
=\sum\limits_{n=1}^\infty 2^{\kappa na} {\bf E} (\xi_n^a)\\
&\le& \sum\limits_{n=1}^\infty 2^{\kappa na} {\bf E} 
\left( \sum\limits_{(x,y)\in D_n^2,|x-y|=2^{-n}}
(\hat \rho (Z_x(\omega), Z_y(\omega))^a \right)  \\
&\le&
\sum\limits_{n=1}^\infty 2^{\kappa na} \cdot d \cdot 2^{dn} \cdot c\cdot 
2^{-n(d+b)}
= cd\sum\limits_{n=1}^\infty 2^{-n(b-a\kappa)}
=\frac{cd2^{a \kappa-b}}{1-2^{a \kappa -b}} < \infty.
\end{eqnarray*}
Hence, there exists $\Omega_0 \in {\mathcal{F}}, \; {\bf P}(\Omega_0)=1$ such that
$$
S(\omega) := \sup_{n\ge 1} (2^{\kappa n} \xi_n(\omega)) < \infty\;\; 
\mbox{ for all } \omega \in \Omega_0.
$$
Further,
$$
{\bf E}(S^a)\le {\bf E}\left(\sum\limits_{n=1}^\infty (2^{\kappa n
} \xi_n)^a  
\right)
\le \frac{cd 2^{a\kappa -b}}{1-2^{a \kappa-b}}.
$$
Let $x, y\in \bigcup\limits_{n=1}^\infty D_n$  such that  $|x-y|_\infty \le r < 
2^{-m}$, 
where $m \in {\bf N}_0$. There exists a sequence 
$$
x=x_1, x_2 \ldots , x_l=y
$$
in $\bigcup\limits_{n=m+1}^\infty D_n$, such that 
for each $i=1, \ldots l -1$ there exists $n(i)\ge m+1$ 
which satisfies $x_i, x_{i+1}\in D_{n(i)}$ and $|x_i-x_{i+1}|=2^{-n(i)}$ 
and 
$$
|\{i \in \{1, \ldots , l-1   \}:n(i)=k  \}| \le 2d\;\; \mbox{ for all } k 
\ge m+1.
$$
For $\omega \in \Omega_0$ and 
$0 < r < 1$ with $2^{-m-1} \le r < 2^{-m}$, we get
\begin{eqnarray*}
&&\sup\{\hat \rho(Z_x(\omega), Z_y(\omega)); x, y\in \bigcup\limits_{n=1}^\infty D_n, 
|x-y|_\infty \le r  \}\\
&\le& 2d \sum\limits_{n=m+1}^\infty \xi_n(\omega) \le 2d 
S(\omega)
\sum\limits_{n=m+1}^\infty 2^{-\kappa n}\\
&=&2^{-\kappa(m+1)} \frac{2d}{1-2^{-\kappa}} S(\omega)
\le \frac{2d}{1-2^{-\kappa}}S(\omega) r^\kappa.
\end{eqnarray*}
The statement in the lemma now follows by the continuity of $Z$. The final statement follows by an 
application of Chebychev's inequality.
\hfill $\Box$\\

\noindent {\bf Proof of Lemma \ref{Chaining}.}  
For each $j \in {\bf N}_0$ and 
each $x\in \Theta_{j+1}$, define  $g_j(x)\in \Theta_j$ such that 
$d(x,g_j(x)) \le \delta_j$ (such a $g_j(x)$ exists due to the 
assumptions in the lemma). We will show, that for each
$x\in \Theta$ there exists a sequence $x_0, x_1, \ldots$ such that
$x=\lim_{j\rightarrow \infty} x_j, x_j\in \Theta_j$ and 
$x_j=g_j(x_{j+1})$ for all $j\in {\bf N}_0$.

To see this, let $\delta_j^*=\sum_{i=j}^\infty \delta_i$ and 
$\tilde{\Theta}_j(x):= \{y \in \Theta_j:d(y,x) \le \delta_j^*  \}$ for  
$x\in \Theta$. Then $\tilde{\Theta}_j(x) \neq \emptyset$ and 
$\tilde{x}\in \tilde{\Theta}_{j+1}(x)$ implies $g_j(\tilde{x})\in \tilde{\Theta}_j(x)$. 
Therefore, there exists a sequence $x_0, x_1, x_2, \ldots $ which satisfies
$x_j\in \tilde{\Theta}_j(x)$ and $x_j=g_j(x_{j+1})$ for all
$j \in {\bf N}_0$. Since 
$\lim_{j\rightarrow \infty} \delta_j^*=0$, we have $x=\lim_{j \to \infty}x_j$. 
We will write  $x_j(x)$ instead of  $x_j$. \\

Fix $x \in \Theta$. The continuity of  $Z$ implies
$$
\hat \rho(Z_x,Z_{x_0}) \le \sum\limits_{j=0}^\infty 
\hat \rho \left(Z_{x_{j+1}(x)},Z_{x_j(x)}\right).
$$
Therefore, 
\begin{eqnarray*}
&\sup\limits_{x\in \Theta} \hat\rho (Z_x, Z_{x_0}) \le \sup\limits_{x \in \Theta}
\sum\limits_{j=0}^\infty \hat \rho \left( Z_{x_{j+1}(x)},Z_{x_j(x)}\right)\\
&\le   \sum\limits_{j=0}^\infty \max\limits_{x_{j+1}\in \Theta_{j+1}} 
\hat \rho \left( Z_{x_{j+1}},Z_{g_j(x_{j+1})} \right).
\end{eqnarray*}
Hence,
\begin{eqnarray*}
{\bf P}\{\sup_{x,y \in \Theta} \hat \rho(Z_x,Z_y) \ge u  \} & \le & 
{\bf P}\{\sup_{x \in \Theta} \hat \rho(Z_x,Z_{x_0}) \ge u/2  \}\\ 
& \le & {\bf P}\{\sum\limits_{j=0}^\infty 
\max\limits_{x\in \Theta_{j+1}}  \hat \rho \left( Z_x,Z_{g_j(x)} \right) 
\ge \frac{u}{2} \sum\limits_{j=0}^\infty \varepsilon_j\}\\
& \le &  
\sum\limits_{j=0}^\infty \sum\limits_{x\in \Theta_{j+1}}{\bf P}\{\hat \rho \left(Z_x,Z_{g_j(x)} \right) \ge 
\varepsilon_j u/2\}\\
& \le & 
\sum\limits_{j=0}^\infty |\Theta_{j+1}| \sup\limits_{d(x,y)\le \delta_j} 
{\bf P}\{\hat \rho \left(Z_x,Z_y\right) \ge \varepsilon_j u/2\}.
\end{eqnarray*}
This completes the proof of the lemma.
\hfill$\Box$\\

\noindent {\bf Proof of Lemma \ref{GRR}.}  The proof is essentially a combination of those of 
\cite{ai} and \cite{dz}. The case $V=0$ is clear (by our conventions about $V/0$), so we assume $V>0$.
We abbreviate
\begin{eqnarray*}
\tilde \Psi (x,y)&:=& \left\{ 
\begin{array}{ll}
\Psi\left( \frac{\hat \rho(f(x),f(y))}{p(d(x,y))} \right),& \mbox{ if } x \neq y,\\
0 & \mbox{ if } x=y.
\end{array}
\right.\\
\end{eqnarray*}
Fix $x \neq y, \,  x,y \in \Theta$ and define $\rho:=d(x,y)$ and
$$ I(u) := \int\limits_{\Theta} \tilde \Psi (u,z) m(dz), \quad u \in \Theta.$$
If either 
$m(K_\varepsilon (x)) = 0$ or $m(K_\varepsilon (y)) = 0$ for some $\varepsilon > 0$ or $V=0$, then
there is nothing to 
show, so we will assume $V > 0$, $m(K_\varepsilon (x)) > 0$ and $m(K_\varepsilon (y)) > 0$ for 
all $\varepsilon > 0$.\\

Let $$U:= \{ z \in \Theta : d(x,z) \leq \rho \mbox{ and } d(y,z) \leq \rho\}.$$
By the definition of $I$ there exists $x_{-1} \in U$ such that
\begin{equation}\label{mu}
I(x_{-1}) \leq \frac{V}{m(U)}.
\end{equation}
Let $\rho=r_{-1} \ge r_0 \ge r_1 \ge ...$ be a sequence of strictly positive reals which we will specify below.
We will recursively define $x_n \in K_{r_n} (x), n \in {\bf N}_0$ such that
\begin{equation}\label{aa}
I (x_n) \leq \frac{2V}{m(K_{r_n} (x))} \qquad \mbox{and}
\end{equation}
\begin{equation}\label{bb}
\tilde \Psi (x_n, x_{n-1}) \leq \frac{2 
I(x_{n-1})}{m(K_{r_n}(x))} , \quad n \in {\bf N}_0. 
\end{equation}
For $n \in {\bf N}_0$ define
$$A_n := \left \{ z \in K_{r_n}(x) : \, I (z) > \frac{2V}{m(K_{r_n}(x))} \quad 
\mbox{or}\quad 
\tilde \Psi  (z,x_{n-1}) > \frac{2 I(x_{n-1})}{m(K_{r_n}(x))} \right \}.$$
Then
$$m(A_n) \leq \frac{m(K_{r_n}(x))}{2V} \int\limits_{\Theta} I(z) m(dz) + 
\frac{m(K_{r_n}(x))}{2I(x_{n-1})} 
\int\limits_{\Theta} \tilde \Psi (z,x_{n-1})m(dz) \leq m(K_{r_n} (x))$$
and the first inequality is strict if $m(A_n) > 0$. In any case we have $m(A_n) 
< m (K_{r_n} (x)).$ 
Now any $x_n \in K_{r_n} (x) \backslash A_n$ \, will satisfy (\ref{aa}) and (\ref{bb}). Using the fact 
that $K_{\rho}(x) \subseteq U$, it follows that
\begin{eqnarray*} 
\hat \rho (f(x_{n+1}), f(x_n)) & \leq & \Psi^{-1} \left( \frac{2 
I(x_n)}{m(K_{r_{n+1}}(x))} \right) p(d(x_n, x_{n+1}))\\
&\leq& \Psi^{-1} \left( \frac{4V}{m(K_{r_{n+1}}(x))m(K_{r_n}(x))}\right)  p(d(x_n, x_{n+1})) \, , n \ge -1.
\end{eqnarray*}
Now, we choose the sequence $r_n$ recursively as follows:
$$
p(2 r_{n+1})=\frac{1}{2}p(r_n+r_{n+1}),\;\;r_{-1}=2\rho.
$$
It is easy to check that this defines the sequence uniquely and that it decreases to zero as $n \to \infty$.
If $n \ge -1$, then
\begin{align*}
&p(d(x_n,x_{n+1})) \le p(d(x_n,x) + d(x,x_{n+1})) \le p(r_n + r_{n+1}) \le 2 p(2r_{n+1}) \\
&\hspace{1cm} = 4 p(2r_{n+1}) - 2  p(2r_{n+1})  \le 4  p(2r_{n+1}) - 4   p(2r_{n+2}).
\end{align*}
Hence,
$$
\hat \rho (f(x_{n+1}), f(x_n)) \le 4 \int_{2r_{n+2}}^{2r_{n+1}} \Psi^{-1} \left( \frac{4V}{(m(K_{s/2}(x)))^2}\right) {\mathrm d} p(s).
$$
The fact that $f$ is continuous (at $x$) implies
\begin{eqnarray*}
\hat \rho (f(x),f(x_{-1})) &\leq& 4 \int\limits^{2r_0}_0 \Psi^{-1} 
\left(\frac{4V}{m(K_{s/2}(x))^2}\right) {\mathrm d} p(s) \\
 &\leq& 4 \int\limits^{4\rho}_0 \Psi^{-1} 
\left(\frac{4V}{m(K_{s/2}(x))^2}\right) {\mathrm d} p(s)    \, .
\end{eqnarray*}
The same estimate holds with $x$ replaced by $y$ (with $y_{-1}:=x_{-1}$). 
Using the triangle inequality we get 
$$
\hat \rho (f(x),f(y)) \leq 8 \max_{z \in \{x,y\}} \int^{4\rho}_0 \Psi^{-1} 
\left(\frac{4V}{m(K_{s/2}(z))^2}\right) {\mathrm d} p(s).   
$$
showing (i) of the Lemma. 
If $N=0$, then there is nothing to show. If $N=1$, then $V \le 1$ and (ii) follows. The general case $N>0$ 
can be reduced to the case $N=1$ by considering the metric $\hat \rho'(x,y):=N^{-1} \hat \rho(x,y)$. \hfill $\Box$


\begin{thebibliography}{0}
%

\bibitem{ai} L.~Arnold and P.~Imkeller, Stratonovich calculus with spatial parameters
and anticipative problems in multiplicative ergodic theory, 
{\it Stoch.~Proc.~Appl.} {\bf 62}, 19--54 (1996).


\bibitem{bh} P.~Baxendale and T.~Harris, Isotropic stochastic flows,
{\it Ann.~Probab.} {\bf 14}, 1155--1179  (1986).

\bibitem{bs} A.~Borodin and P.~Salminen, {\it Handbook of Brownian Motion -- Facts and Formulae}, 
Birkh\"auser, Basel, 1996.


\bibitem{cs} M.~Cranston and M.~Scheutzow, Dispersion rates under finite mode Kolmogorov
flows, {\it Ann.~Appl.~Probab.}, {\bf 12}, 511--532 (2002).

\bibitem{css1} M.~Cranston, M.~Scheutzow and D.~Steinsaltz, Linear expansion of isotropic
Brownian flows, {\it Electron.~Commun.~Probab.} {\bf 4}, 91--101 (1999).

\bibitem{css2} M.~Cranston, M.~Scheutzow and D.~Steinsaltz, Linear bounds for stochastic
dispersion, {\it Ann.~Probab.} {\bf 28}, 1852--1869 (2000).

\bibitem{dkn} R.~Dalang, D.~Khoshnevisan, and E.~Nualart, 
Hitting probabilities for systems of non-linear stochastic heat equations with additive noise,
{\it Alea} {\bf 3}, 231--271 (2007).

\bibitem{dz} G.~Da Prato and J.~Zabczyk, {\it Ergodicity for Infinite Dimensional Systems},
Cambridge University Press, 1996.

\bibitem{dembo} A.~Dembo and O.~Zeitouni, {\it Large Deviations Techniques and Applications},
2nd edition, Springer, New York, 1998.

\bibitem{ds}  G.~Dimitroff and M.~Scheutzow, Attractors and expansion for 
Brownian flows, submitted.

\bibitem{f} W.~Feller, The asymptotic distribution of the range of sums of independent random variables, 
{\it Ann.~Math.~Statistics} {\bf 22}. 427--432 (1951).

\bibitem{grr} A.M.~Garsia, E.~Rodemich, and H.~Rumsey, A real variable lemma and the continuity 
of paths of some Gaussian processes, 
{\it Indiana Univ.~Math.~Journal} {\bf 20}, 565--578 (1970).

\bibitem{kal} O.~Kallenberg, {\it Foundations of Modern Probability}, 
2nd edition, Springer, 2002.

\bibitem{k} H.~Kunita, {\it Stochastic Flows and Stochastic Differential Equations},
Cambridge University Press, 1990.

\bibitem{lt} M.~Ledoux and M.~Talagrand, {\it Probability in Banach Spaces}, Springer, 1991.

\bibitem{ls01} H.~Lisei and M.~Scheutzow, Linear bounds and Gaussian tails in a
stochastic dispersion model, {\it Stochastics and Dynamics} {\bf 1}, 389--403 (2001).

\bibitem{ls03} H.~Lisei and M.~Scheutzow, {\it On the dispersion of sets under the action of an 
isotropic Brownian flow},
in: Proceedings of the Swansea 2002 Workshop {\it Probabilistic Methods in Fluids},
ed.~I.~Davies, 224--238, World Scientific, 2003.

\bibitem{p} D.~Pollard, {\it Empirical Processes: Theory and Applications}, IMS, 1990.

\bibitem{s07} M.~Scheutzow, Attractors for Ergodic and Monotone Random Dynamical Systems,
in:  {\it Seminar on Stochastic Analysis, Random Fields and
Applications V},  ed.~R.~Dalang, M.~Dozzi, F.~Russo, 331--344, Birkh\"auser, 2007. 

\bibitem{ss} M.~Scheutzow and D.~Steinsaltz, Chasing balls through martingale fields,
{\it Ann.~Probab.} {\bf 30}, 2046--2080 (2002).

\bibitem{T05} M.~Talagrand, {\it The Generic Chaining}, Springer, 2005.

\bibitem{walsh} J.~Walsh, An introduction to stochastic partial differential equations, in: {\it \'Ecole d'\'et\'e 
de probabilit\'es de Saint-Flour XIV - 1984}, Lect. Notes Math. 1180, 265--437, Springer, 1986.

\end{thebibliography}
\end{document}